# Two Remarkable Ortho-Homological Triangles

Prof. Ion Pătraşcu - The National College "Fraţii Buzeşti", Craiova, Romania
Prof. Florentin Smarandache – University of New Mexico, U.S.A.

In a previous paper [5] we have introduced the ortho-homological triangles, which are triangles that are orthological and homological simultaneously.

In this article we call attention to two remarkable ortho-homological triangles (the given triangle $ABC$ and its first Brocard's triangle), and using the Sondat's theorem relative to orthological triangles, we emphasize on four important collinear points in the geometry of the triangle. Orthological / homological / orthohomological triangles in the 2D-space are generalized to orthological / homological / orthohomological polygons in 2D-space, and even more to orthological / homological / orthohomological triangles, polygons, and polyhedrons in 3D-space.

**Definition 1**

The first Brocard triangle of a given triangle $ABC$ is the triangle formed by the projections of the symmedian center of the triangle $ABC$ on its perpendicular bisectors.

**Observation**

In figure 1 we note with $K$ the symmedian center, $OA'$, $OB'$, $OC'$ the perpendicular bisectors of the triangle $ABC$ and $A_1B_1C_1$ the first Brocard's triangle.

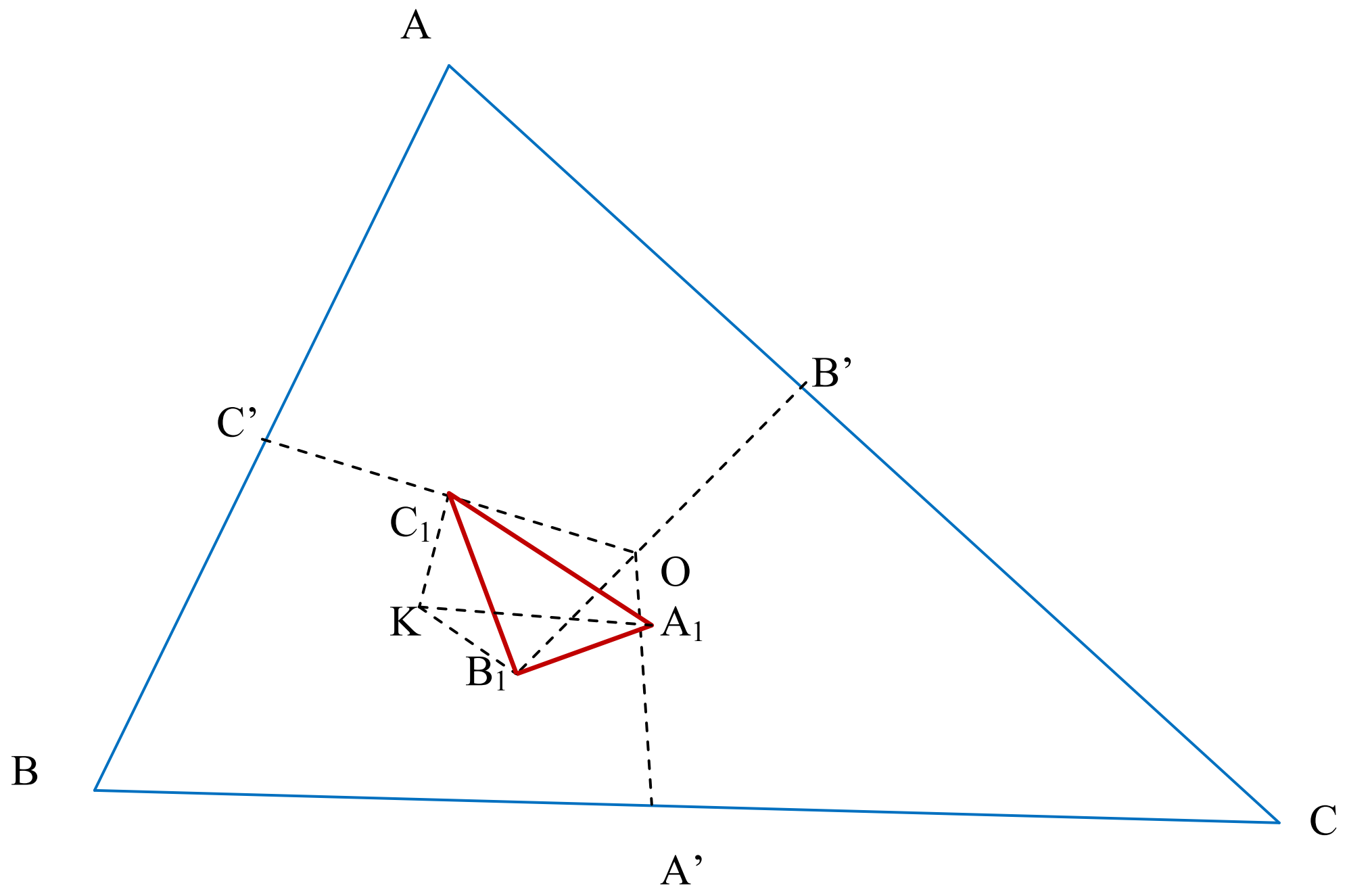

Fig. 1

**Theorem 1**

If $ABC$ is a given triangle and $A_1B_1C_1$ is its first triangle Brocard, then the triangles $ABC$ and $A_1B_1C_1$ are ortho-homological.

We'll perform the proof of this theorem in two stages.

I. We prove that the triangles $A_1B_1C_1$ and $ABC$ are orthological.

The perpendiculars from $A_1$, $B_1$, $C_1$ on $BC$, $CA$ respective $AB$ are perpendicular bisectors in the triangle $ABC$, therefore are concurrent in $O$, the center of the circumscribed circle of triangle $ABC$ which is the orthological center for triangles $A_1B_1C_1$ and $ABC$.

II. We prove that the triangles $A_1B_1C_1$ and $ABC$ are homological, that is the lines $AA_1$, $BB_1$, $CC_1$ are concurrent.

To continue with these proves we need to refresh some knowledge and some helpful results.

**Definition 2**

In any triangle $ABC$ there exist the points $\Omega$ and $\Omega'$ and the angle $\omega$ such that:

$$m(\sphericalangle \Omega AB) = \sphericalangle \Omega BC = \sphericalangle \Omega CA = \omega$$

$$m(\sphericalangle \Omega' BA) = \sphericalangle \Omega' CA = \sphericalangle \Omega' AB = \omega$$

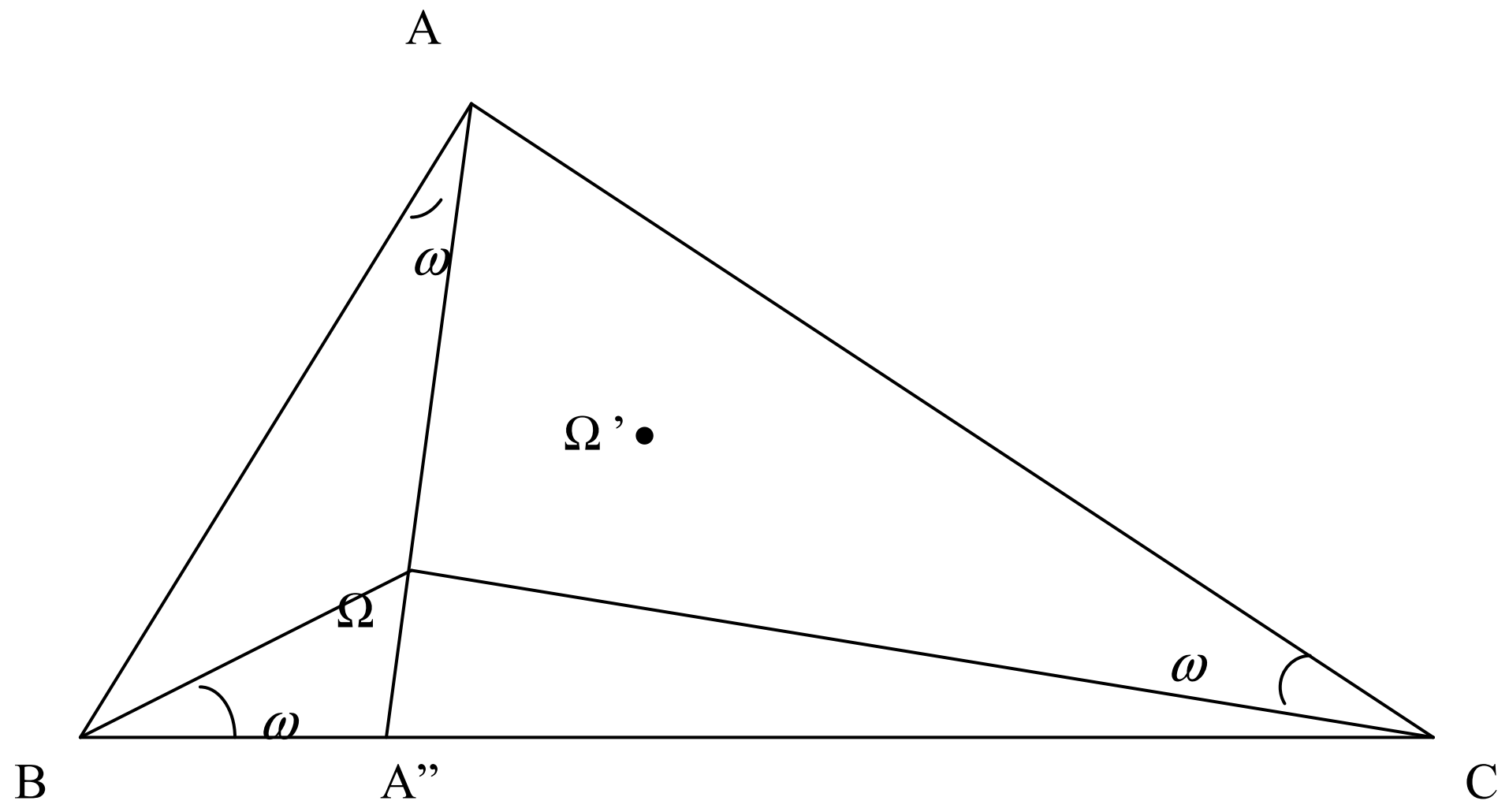

Fig. 2

The points $\Omega$ and $\Omega'$ are called the first, respectively the second point of Brocard and $\omega$ is called the Brocard's angle.

**Lemma 1**

In the triangle $ABC$ let $\Omega$ the first point of Brocard and $\{A''\}$ $\{A''\} = A\Omega \cap BC$, then:

$$\frac{BA"}{CA"} = \frac{c^2}{a^2}$$

**Proof**

$$Aria_{\triangle}ABA" = \frac{1}{2} AB \cdot AA" \sin \omega \quad (1)$$

$$Aria_{\triangle}ACA" = \frac{1}{2} AC \cdot AA" \sin(A - \omega) \quad (2)$$

From (1) and (2) we find:

$$\frac{Aria_{\triangle}ABA"}{Aria_{\triangle}ACA"} = \frac{AB \cdot \sin \omega}{AC \cdot \sin(A-\omega)} \quad (3)$$

On the other side, the mentioned triangles have the same height built from $A$, therefore:

$$\frac{Aria_{\triangle}ABA"}{Aria_{\triangle}ACA"} = \frac{BA"}{CA"} \quad (4)$$

From (3) and (4) we have:

$$\frac{BA"}{CA"} = \frac{AB \cdot \sin \omega}{AC \cdot \sin(A-\omega)} \quad (5)$$

Applying the sinus theorem in the triangle $A\Omega C$ and in the triangle $B\Omega C$, it results:

$$\frac{C\Omega}{\sin(A-\omega)} = \frac{AC}{\sin A\Omega C} \quad (6)$$

$$\frac{C\Omega}{\sin \omega} = \frac{BC}{\sin B\Omega C} \quad (7)$$

Because

$$m(\sphericalangle A\Omega C) = 180^o - A$$

$$m(\sphericalangle B\Omega C) = 180^o - C$$

From the relations (6) and (7) we find:

$$\frac{\sin \omega}{\sin(A-\omega)} = \frac{AC}{BC} \cdot \frac{\sin C}{\sin A} \quad (8)$$

Applying the sinus theorem in the triangle $ABC$ leads to:

$$\frac{\sin C}{\sin A} = \frac{AB}{BC} \quad (9)$$

The relations (5), (8), (9) provide us the relation:

$$\frac{BA"}{CA"} = \frac{c^2}{a^2}$$

**Remark 1**

By making the notations: $\{B"\} = B\Omega C \cap AC$ and $\{C"\} = C\Omega A \cap AB$ we obtain also the relations:

$$\frac{CB"}{AB"} = \frac{a^2}{b^2} \text{ and } \frac{AC"}{BC"} = \frac{b^2}{c^2}$$

**Lemma 2**

In a triangle $ABC$, the Brocard's Cevian $B\Omega$, symmedian from $C$ and the median from $A$ are concurrent.

**Proof**

It is known that the symmedian $CK$ of triangle $ABC$ intersects $AB$ in the point $C_2$ such that $\frac{AC_2}{BC_2}=\frac{b^2}{c^2}$. We had that the Cevian $B\Omega$ intersects $AC$ in $B''$ such that $\frac{BC''}{B''A}=\frac{a^2}{b^2}$.

The median from $A$ intersects $BC$ in $A'$ and $BA'=CA'$.

Because $\frac{A'B}{A'C}\cdot\frac{B''C}{B''A}\cdot\frac{C_2A}{C_2B}=1$, the reciprocal of Ceva's theorem ensures the concurrency of the lines $B\Omega,\ CK$ and $AA'$.

**Lemma 3**

Give a triangle $ABC$ and $\omega$ the Brocard's angle, then

$$ctg\,\omega = ctgA + ctgB + ctgC \tag{9}$$

**Proof**

From the relation (8) we find:

$$\sin(A-\omega)=\frac{a}{b}\cdot\frac{\sin A}{\sin C}\cdot\sin\omega \tag{10}$$

From the sinus' theorem in the triangle $ABC$ we have that

$$\frac{a}{b}=\frac{\sin A}{\sin B}$$

Substituting it in (10) it results: $\sin(A-\omega)=\frac{\sin^2 A\cdot\sin\omega}{\sin B\cdot\sin C}$

Furthermore we have:

$$\sin(A-\omega)=\sin A\cdot\cos\omega-\sin\omega\cdot\cos A$$

$$\sin A\cdot\cos\omega-\sin\omega\cdot\cos A=\frac{\sin^2 A\cdot\sin\omega}{\sin B\cdot\sin C} \tag{11}$$

Dividing relation (11) by $\sin A\cdot\sin\omega$ and taking into account that $\sin A=\sin(B+C)$, and $\sin(B+C)=\sin B\cdot\cos C+\sin C\cdot\cos B$ we obtain relation (5)

**Lemma 4**

If in the triangle $ABC$, $K$ is the symmedian center and $K_1,K_2,K_3$ are its projections on the sides $BC,\ CA,\ AB$, then:

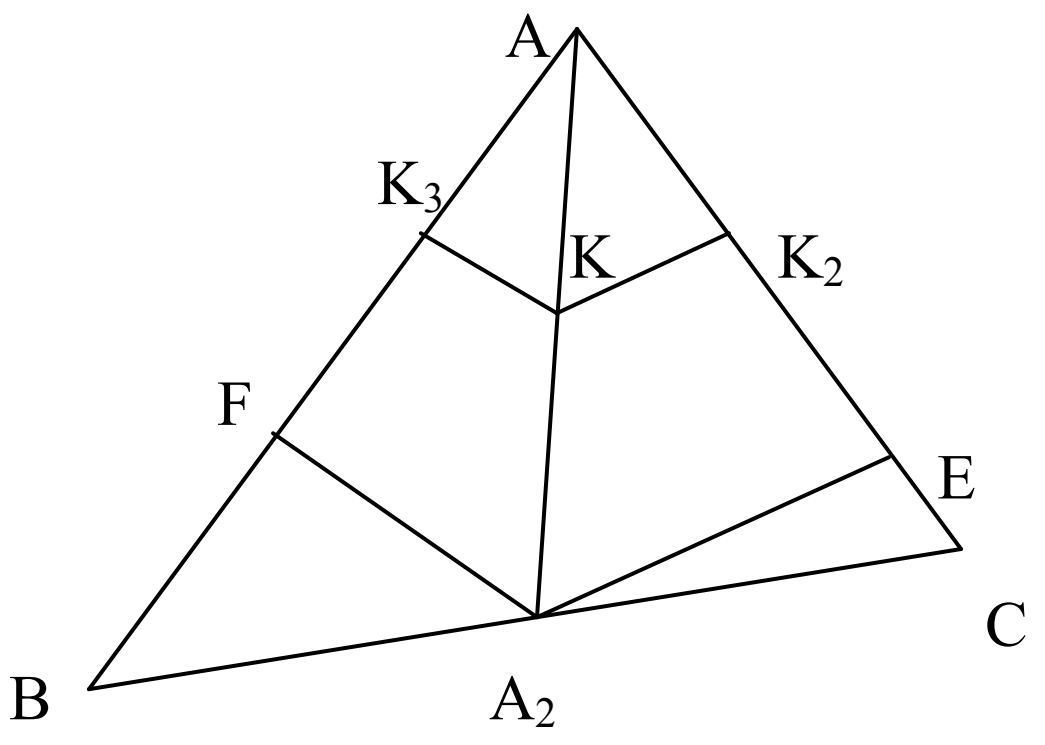

Fig. 3

$$\frac{KK_1}{a}=\frac{KK_2}{b}=\frac{KK_3}{c}=\frac{1}{2}tg\omega$$

**Proof:**

Let $AA_2$ the symmedian in the triangle $ABC$, we have:

$$\frac{BA_2}{CA_2}=\frac{Aria\triangle BAA_2}{Aria\triangle CAA_2},$$

where $E$ and $F$ are the projection of $A_2$ on $AC$ respectively $AB$.

It results that $\frac{A_2F}{A_2E}=\frac{c}{b}$

From the fact that $\triangle AKK_3 \sim \triangle AA_2F$ and $\triangle AKK_2 \sim \triangle AA_2E$ we find that $\frac{KK_3}{KK_2}=\frac{A_2F}{A_2E}$

Also: $\frac{KK_2}{b}=\frac{KK_3}{c}$, and similarly: $\frac{KK_1}{a}=\frac{KK_2}{b}$, consequently:

$$\frac{KK_1}{a}=\frac{KK_2}{b}=\frac{KK_3}{c} \qquad (12)$$

The relation (12) is equivalent to:

$$\frac{aKK_1}{a^2}=\frac{bKK_2}{b^2}=\frac{cKK_3}{c^2}=\frac{aKK_1+bKK_2+cKK_3}{a^2+b^2+c^2}$$

Because

$$aKK_1+bKK_2+cKK_3=2Aria\triangle ABC=2S,$$

we have:

$$\frac{KK_1}{a}=\frac{KK_2}{b}=\frac{KK_3}{c}=\frac{2S}{a^2+b^2+c^2}$$

If we note $H_1, H_2, H_3$ the projections of $A, B, C$ on $BC, CA, AB$, we have

$$ctgA=\frac{H_2A}{BH_2}=\frac{bc\cos A}{2S}$$

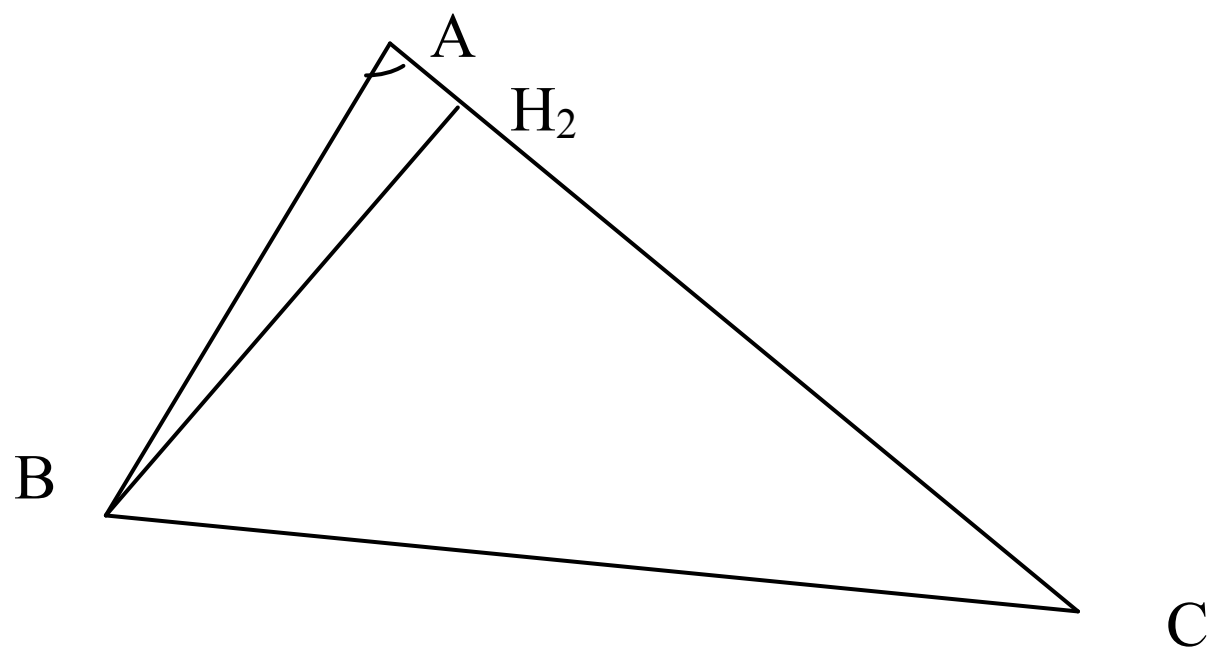

Fig. 4

From the cosine's theorem it results that : $b \cdot c \cdot \cos A = \frac{b^2 + c^2 - a^2}{2}$, and therefore

$$ctgA = \frac{b^2 + c^2 - a^2}{4S}$$

Taking into account the relation (9), we find:

$$ctg\omega = \frac{a^2 + b^2 + c^2}{4S},$$

then

$$tg\omega = \frac{4S}{a^2 + b^2 + c^2}$$

and then

$$\frac{KK_1}{a} = \frac{2S}{a^2 + b^2 + c^2} = \frac{1}{2} tg\omega .$$

**Lemma 5**

The Cevians $AA_1$, $BB_1$, $CC_1$ are the isotomics of the symmedians $AA_2$, $BB_2$, $CC_2$ in the triangle $ABC$.

**Proof:**

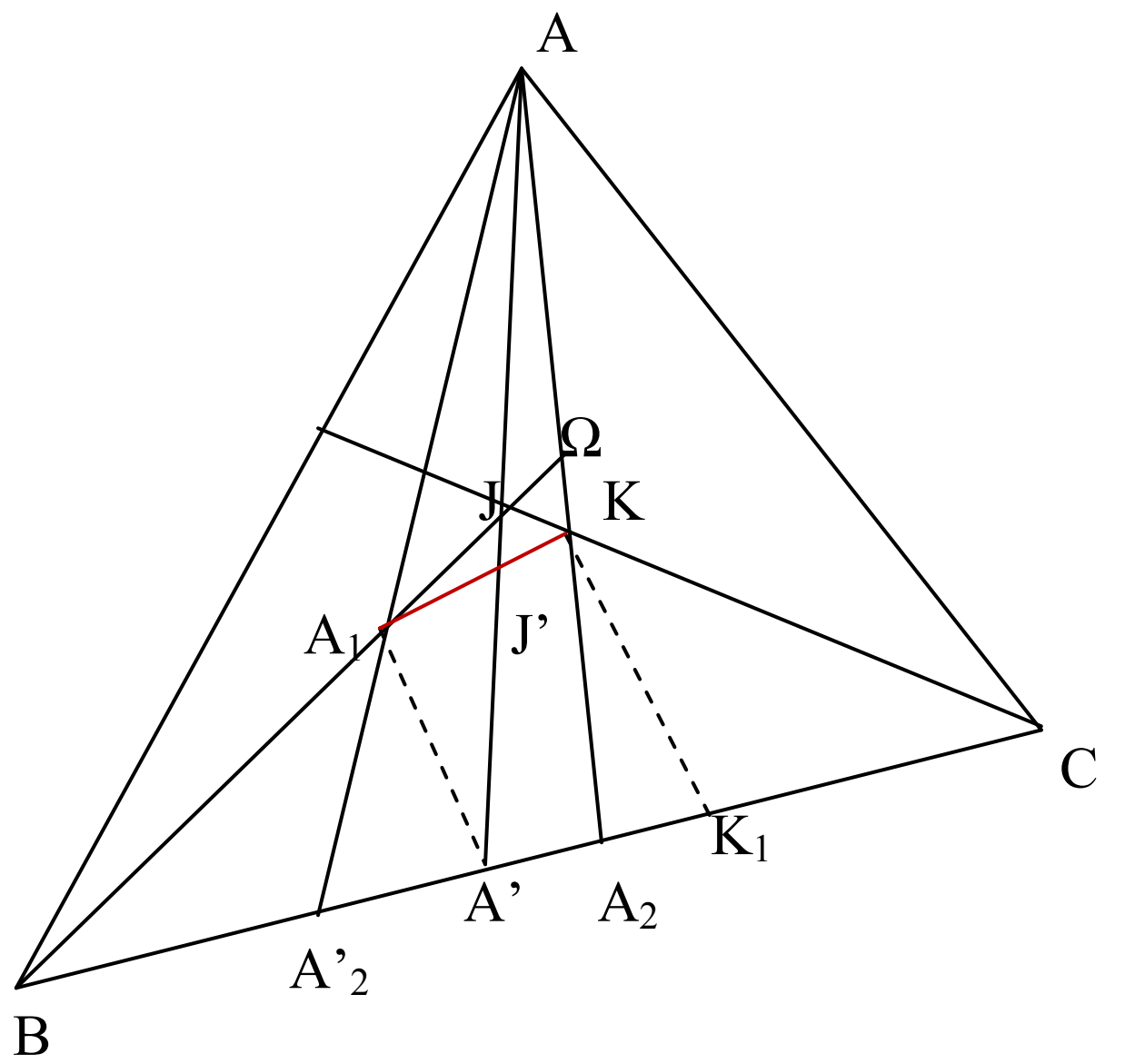

Fig. 5

In figure 5 we note $J$ the intersection point of the Cevians from the Lemma 2.

Because $KA_1 \parallel BC$, we have that $A_1A' = KK_1 = \frac{1}{2} a\,tg\,\omega$. On the other side from the right triangle $A'A_1B$ we have: $tg\sphericalangle A_1BA' = \frac{A_1A'}{BA'} = tg\,\omega$, consequently the point $A_1$, the vertex of the first triangle of Brocard belongs to the Cevians $B\Omega$.

We note $\{J'\} = A_1K \cap AA'$, and evidently from $A_1K \parallel BC$ it results that $JJ'$ is the median in the triangle $JA_1K$, therefore $A_1J' = J'K$.

We note with $A'_2$ the intersection of the Cevians $AA_1$ with $BC$, because $A_1K \parallel A'_2A_2$ and $AJ'$ is a median in the triangle $AA_1K$ it results that $AA'$ is a median in triangle $AA'_2A_2$ therefore the points $A'_2$ and $A_2$ are isometric.

Similarly it can be shown that $BB_2'$ and $CC_2'$ are the isometrics of the symmedians $BB_2$ and $CC_2$.

The second part of this proof: Indeed it is known that the isometric Cevians of certain concurrent Cevians are concurrent and from Lemma 5 along with the fact that the symmedians of a triangle are concurrent, it results the concurrency of the Cevians $AA_1$, $BB_1$, $CC_1$ and therefore the triangle $ABC$ and the first triangle of Brocard are homological. The homology's center (the concurrency point) of these Cevians is marked in some works with $\Omega''$ with and it is called the third point of Brocard.

From the previous proof, it results that $\Omega''$ is the isotomic conjugate of the symmedian center $K$.

**Remark 2**

The triangles ABC and $A_1B_1C_1$ (first Brocard triangle) are *triple-homological*, since first time the Cevians $AB_1$, $BC_1$, $CA_1$ are concurrent (in a Brocard point), second time the Cevians $AC_1$, $BA_1$, $CB_1$ are also concurrent (in the second Brocard point), and third time the Cevians $AA_1$, $BB_1$, $CC_1$ are concurrent as well (in the third point of Brocard).

**Definition 3**

It is called the Tarry point of a triangle $ABC$, the concurrency point of the perpendiculars from $A$, $B$, $C$ on the sides $B_1C_1$, $C_1A_1$, $A_1B_1$ of the Brocard's first triangle.

**Remark 3**

The fact that the perpendiculars from the above definition are concurrent results from the theorem 1 and from the theorem that states that the relation of triangles' orthology is symmetric.

We continue to prove the concurrency using another approach that will introduce supplementary information about the Tarry's point.

We'll use the following:

**Lemma 6:**

The first triangle Brocard of a triangle and the triangle itself are similar.

**Proof**

From $KA_1 \parallel BC$ and $OA' \perp BC$ it results that

$$m(\sphericalangle KA_1O) = 90°$$

(see Fig. 1), similarly

$$m(\sphericalangle KB_1O) = m(\sphericalangle KC_1O) = 90°$$

and therefore the first triangle of Brocard is inscribed in the circle with $OK$ as diameter (this circle is called the Brocard circle).

Because

$$m(\sphericalangle A_1OC_1) = 180° - B$$

and $A_1, B_1, C_1, O$ are concyclic, it results that $\sphericalangle A_1B_1C_1 = \sphericalangle B$, similarly

$$m(\sphericalangle B'OC') = 180° - A,$$

it results that

$$m(\sphericalangle B_1OC_1) = m(A)$$

but

$$\sphericalangle B_1OC_1 \equiv \sphericalangle B_1A_1C_1,$$

therefore

$$\sphericalangle B_1A_1C_1 = \sphericalangle A$$

and the triangle $A_1B_1C_1$ is similar wit the triangle $ABC$.

**Theorem 2**

The orthology center of the triangle $ABC$ and of the first triangle of Brocard is the Tarry's point $T$ of the triangle $ABC$, and $T$ belongs to the circumscribed circle of the triangle $ABC$.

**Proof**

We mark with $T$ the intersection of the perpendicular raised from $B$ on $A_1C_1$ with the perpendicular raised from $C$ on $A_1B_1$ and let

$$\{B_1'\} = BT \cap A_1C_1,\ A_1\{C_1'\} = A_1B_1 \cap CT.$$

We have

$$m(\triangle B_1'TC_1') = 180° - m(\triangle C_1A_1B_1)$$

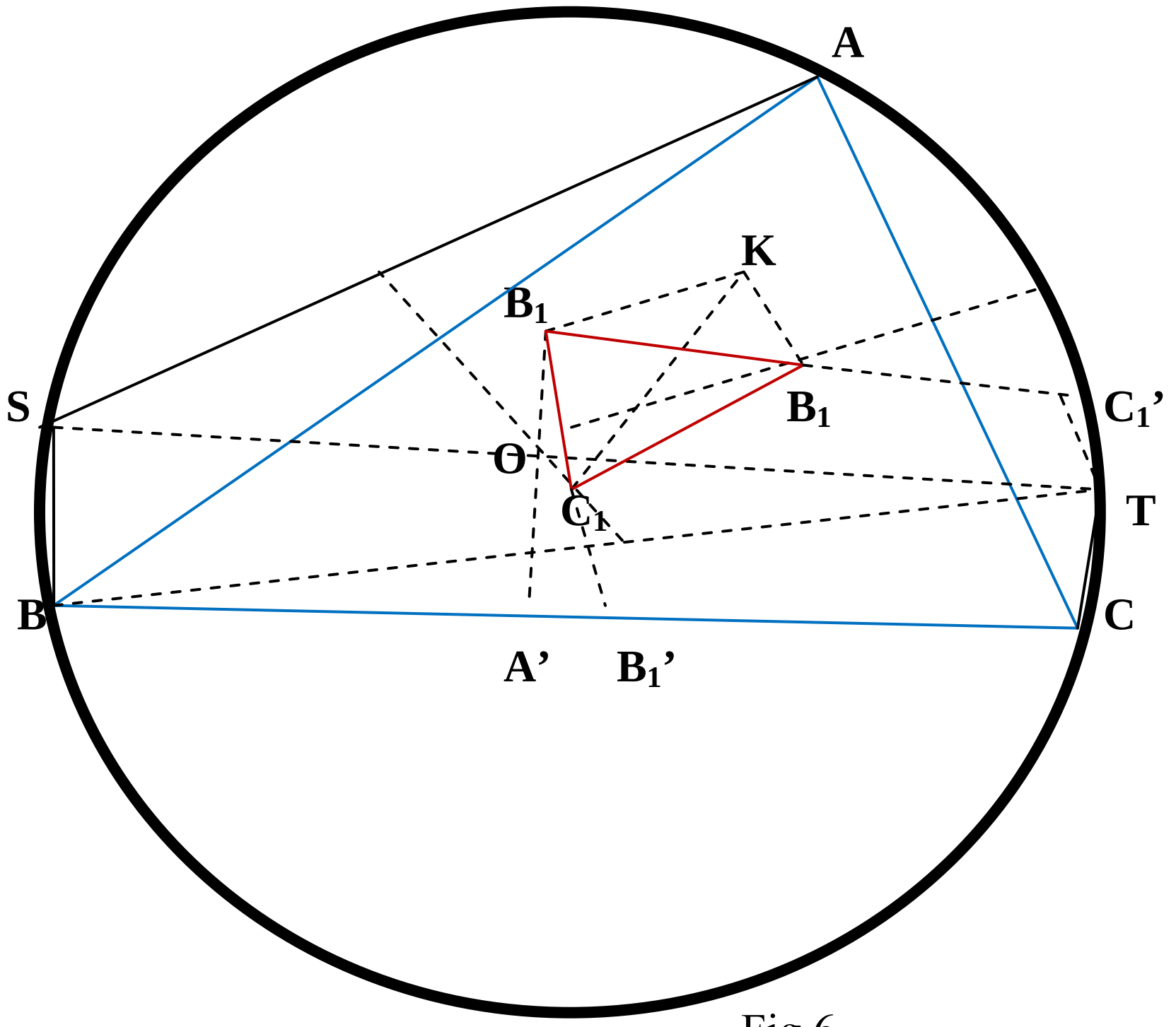

Fig 6

But because of Lemma 6 $\sphericalangle C_1A_1B_1 = \sphericalangle A$.

It results that $m(\sphericalangle B_1'TC_1') = 180° - A$, therefore

$$m(\sphericalangle BTC') + m(\sphericalangle BAC) = 180°$$

Therefore $T$ belongs to the circumscribed circle of triangle $ABC$

If $\{A_1'\} = B_1C_1 \cap AT$ and if we note with $T'$ the intersection of the perpendicular raised from $A$ on $B_1C_1$ with the perpendicular raised from $B$ on $A_1C_1$, we observe that

$$m(\sphericalangle B_1'T'A_1') = m(\sphericalangle A_1C_1B_1)$$

therefore

$$m(\sphericalangle BT'A) + m(\sphericalangle BCA)$$

and it results that $T'$ belongs to the circumscribed triangle $ABC$.

Therefore $T = T'$ and the theorem is proved.

**Theorem 3**

If through the vertexes $A, B, C$ of a triangle are constructed the parallels to the sides $B_1C_1, C_1A_1$ respectively $A_1B_1$ of the first triangle of Brocard of this triangle, then these lines are concurrent in a point S (the Steiner point of the triangle)

**Proof**

We note with $S$ the polar intersection constructed through $A$ to $B_1C_1$ with the polar constructed through $B$ to $A_1C_1$ (see Fig. 6).

We have

$$m(\sphericalangle ASB) = 180° - m(\sphericalangle A_1C_1B_1) \text{ (angles with parallel sides)}$$

because

$$m(\sphericalangle A_1C_1B_1) = m\sphericalangle C,$$

we have

$$m(\sphericalangle ASB) = 180° - m\sphericalangle C,$$

therefore $A_1SB_1C$ are concyclic.

Similarly, if we note with $S'$ the intersection of the polar constructed through $A$ to $B_1C_1$ with the parallel constructed through $C$ to $A_1B_1$ we find that the points $A_1S_1'B_1C$ are concyclic.

Because the parallels from $A$ to $B_1C_1$ contain the points $A, S, S'$ and the points $S, S', A$ are on the circumscribed circle of the triangle, it results that $S = S'$ and the theorem is proved.

**Remark 4**

Because $SA \| B_1C_1$ and $B_1C_1 \perp AT$, it results that

$$m(\sphericalangle SAT) = 90°,$$

but $S$ and $T$ belong to the circumscribed circle to the triangle $ABC$, consequently the Steiner's point and the Tarry point are diametric opposed.

**Theorem 4**

In a triangle $ABC$ the Tarry point $T$, the center of the circumscribed circle $O$, the third point of Brocard $\Omega''$ and Steiner's point $S$ are collinear points

**Proof**

The P. Sondat's theorem relative to the orthological triangles (see [4]) says that the points $T$, $O$, $\Omega''$ are collinear, therefore the points: $T$, $O$, $\Omega''$, $S$ are collinear.

**Open Questions**

1) Is it possible to have two triangles which are four times, five times, or even six times orthological? But triangles which are four times, five times, or even six times homological? What about orthohomological? What is the largest such rank?
   For two triangles $A_1B_1C_1$ and $A_2B_2C_2$, we can have (in the case of orthology, and similarity in the cases of homology and orthohomology) the following 6 possibilities:
   1) the perpendicular from $A_1$ onto $B_2C_2$, the perpendicular from $B_1$ onto $C_2A_2$, and the perpendicular from $C_1$ onto $A_2B_2$ concurrent;

2) the perpendicular from $A_1$ onto $B_2C_2$, the perpendicular from $B_1$ onto $A_2B_2$, and the perpendicular from $C_1$ onto $C_2A_2$ concurrent;
3) the perpendicular from $B_1$ onto $B_2C_2$, the perpendicular from $A_1$ onto $C_2A_2$, and the perpendicular from $C_1$ onto $A_2B_2$ concurrent;
4) the perpendicular from $B_1$ onto $B_2C_2$, the perpendicular from $A_1$ onto $A_2B_2$, and the perpendicular from $C_1$ onto $C_2A_2$ concurrent;
5) the perpendicular from $C_1$ onto $B_2C_2$, the perpendicular from $B_1$ onto $C_2A_2$, and the perpendicular from $A_1$ onto $A_2B_2$ concurrent;
6) the perpendicular from $C_1$ onto $B_2C_2$, the perpendicular from $B_1$ onto $A_2B_2$, and the perpendicular from $A_1$ onto $C_2A_2$ concurrent.

2) We generalize the orthological, homological, and orthohomological triangles to respectively orthological, homological, and orthohomological polygons and polyhedrons. Can we have double, triple, etc. orthological, homological, or orthohomological polygons and polyhedrons? What would be the largest rank for each case?

3) Let's have two triangles in a plane. Is it possible by changing their positions in the plane and to have these triangles be orthological, homological, orthohomological? What is the largest rank they may have in each case?

4) Study the orthology, homology, orthohomology of triangles and poligons in a 3D space.

5) Let's have two triangles, respectively two polygons, in a 3D space. Is it possible by changing their positions in the 3D space to have these triangles, respectively polygons, be orthological, homological, or orthohological?
Similar question for two polyhedrons?